\documentclass{article}
\usepackage[ansinew]{inputenc}
\usepackage[english]{babel}
\usepackage{latexsym}
\usepackage{amsthm}
\usepackage{amsmath}
\usepackage{amsfonts}
\usepackage{amssymb}
\usepackage{graphicx}
\usepackage{times}
\usepackage{url}
\usepackage{float}
\usepackage{epic}

\begin{document}

%%%macros
\newcommand{\m}{\medskip}

\newcommand{\CC}{\mathcal{C}}
\newcommand{\OO}{\mathcal{O}}

\newcommand{\dd}{\operatorname{d}}
\newcommand{\e}{\operatorname{e}}
\newcommand{\rr}{\operatorname{r}}

\newcommand{\A}{\mathbb{A}}
\renewcommand{\P}{\mathbb{P}}
\newcommand{\T}{\mathbb{T}}
\newcommand{\R}{\mathbb{R}}
\newcommand{\N}{\mathbb{N}}
\newcommand{\Z}{\mathbb{Z}}

\newcommand{\realamp}{\R\cup\{-\infty\}}

\newcommand{\ptop}[3]{[#1,#2,#3]}
\newcommand{\pto}[2]{{(#1,#2)}}

\newcommand{\im}{\operatorname{im}}
\newcommand{\conv}{\operatorname{conv}}

%%%AQUI VOY PONIENDO LOS MACROS REVISADOS

%%matriz diagonal asociada a la letra a
\newcommand{\diag}[1]{\left(\begin{array}{ccc}
#1 _{11}/2\\
-\infty&#1 _{22}/2\\
-\infty&-\infty&#1 _{33}/2\\
\end{array}\right)}

\newcommand{\SHP}{\operatorname{shape}}

\newcommand{\suces}[3]{{#1_{#2},\hdots,#1_{#3}}}

\newcommand{\sucespar}[1]{#1_{21},#1_{32},#1_{31}}
\newcommand{\sucesimp}[1]{#1_{1},#1_{2},#1_{3}}

%puntos ASOCIADOS a matriz
\newcommand{\hh}[1]{[-#1 _{11}, -#1 _{21}, -#1 _{31}]}
\newcommand{\vv}[1]{[-#1 _{21}, -#1 _{22}, -#1 _{32}]}
\newcommand{\ii}[1]{[-#1 _{31}, -#1 _{32}, -#1 _{33}]}
\newcommand{\cc}[1]{[#1 _{32}, #1 _{31}, #1 _{21}]}

%% polinomio tropical  de grado dos
\newcommand{\poltrohdos}[1]{#1_{11}\odot X^{\odot2}\oplus #1_{22}\odot
Y^{\odot2}\oplus #1_{33}\odot Z^{\odot2}\oplus #1 _{21}\odot X
\odot Y\oplus #1 _{32}\odot Y \odot Z\oplus #1 _{31}\odot X
\odot Z}

\newcommand{\polly}[1]{#1 {11}\odot X^{\odot2}\oplus #1 {22}\odot
Y^{\odot2}\oplus #1 {33}\odot Z^{\odot2}}

%% polinomio tropical  de grado dos, version maximo
\newcommand{\maxxi}[1]{\max\left\{
#1_{11}+2X, #1_{22}+2Y, #1_{33}+2Z, #1_{21}+X+Y , #1_{32}+Y+Z, #1_{31}+X+Z
\right\}}

\newcommand{\polshp}[1]{
X^{\odot2}\oplus Y^{\odot2}\oplus Z^{\odot2}
 \oplus
{#1 _{21}}\odot X \odot Y\oplus{#1 _{32}}
\odot Y \odot Z \oplus {#1 _{31}}\odot X\odot Z}

\newcommand{\polshpp}[1]{
X^{\odot2}\oplus Y^{\odot2}\oplus Z^{\odot2}
 \oplus
 X \odot Y\oplus {#1 _{32}^+}\odot Y \odot
Z \oplus {#1 _{31}^+}\odot X\odot Z}

\newcommand{\inimp}{\in\{1,2,3\}}

%%matriz asociada a la letra a
\newcommand{\mat}[1]{
\left(\begin{array}{ccc}
#1 _{11}\\
#1 _{21}&#1 _{22}\\
#1 _{31}&#1 _{32}&#1 _{33}\\
\end{array}\right)}

%%%%AQUI ACABAN LOS MACROS REVISADOS

%matriz de forma (shape) asociada a la letra f
\newcommand{\shp}[1]{
\dfrac{1}{2}\left(\begin{array}{ccc}
0\\
#1 _2&0\\
#1 _6&#1 _4&0\\
\end{array}\right)}

\newcommand{\inpar}{\in\{2,4,6\}}

%% polinomio tropical  de grado dos, version maximo, con ceros en diagonal
\newcommand{\maxx}[1]{\max\left\{
2X, 2Y, 2Z, X+Y +\dfrac{#1_{21}}{2}, Y+Z +\dfrac{#1_{32}}{2},
X+Z +\dfrac{#1_{31}}{2}
 \right\}}

\newtheorem{thm}{Theorem}
\newtheorem{prop}{Proposition}
\newtheorem{lem}{Lemma}
\newtheorem{dfn}{Definition}
\newtheorem{cor}{Corollary}
\newtheorem{ex}{Example}
\newtheorem{rem}{Remark}

%%%macros

\title{Tropical conics for the layman}
\author{M. Ansola, M.J. de la
Puente\thanks{Departamento de Algebra, Facultad de Matem{\'a}ticas,
Universidad Complutense, 28040--Madrid, Spain}\ \thanks{Partially
supported by MTM 2005--02865 and UCM 910444.}}
\maketitle
\date{}
\begin{abstract} We present a  simple and elementary procedure to
sketch the tropical conic given by a degree--two
homogeneous tropical polynomial. These conics are trees of
a very particular kind. Given such a tree, we explain how
to compute a defining polynomial.
 Finally, we characterize
those degree--two tropical polynomials which are reducible
and factorize them. We show that there exist irreducible
degree--two tropical polynomials giving rise to  pairs of
tropical lines.
\end{abstract}

%%comments for arXiv
%19 pages, 4 figures.
%The paper has been fully reorganized and the presentation has been improved.
%The title of  former versions of this paper is
%"Metric invariants of tropical conics and factorization of degree--two homogeneous tropical polynomials in three variables".
%To appear in Idempotent and tropical mathematics and problems of mathematical physics (vol.  II), G. Litvinov, V. Maslov, S. Sergeev (eds.), Proceedings of an international Workshop held at the Independent U. Moscow, Russia, 2007.

%There is a figure in a separate file.

2000 Math. Subj. Class.: 05C05;  12K99


\section{Introduction}

In recent years, there has been a growing interest in
\emph{projective tropical geometry},
\cite{Ansola_tri, Develin, Einsiedler, Gathmann,  Itenberg, Jensen, Joswig, Litvinov_ed, Mikhalkin_E, Mikhalkin_J,
Mikhalkin_T, Mikhalkin_W, Richter, Shustin, Shustin_W, Sturmfels,
Viro_D, Viro_W}. This new geometry is related to  toric geometry, \cite{Gelfand, Hsie, Nishinou}.
Several authors have searched for tropical
versions of some classical theorems of projective geometry,
 \cite{Gathmann_M, Speyer, Tabera_Pap, Tabera_Trans, Vigeland_G, Vigeland_S}. Some of these old
theorems involve conics.

The aim of this paper is to present tropical conics to
non--experts, using only tropical algebra (also called max--algebra, max--plus algebra, semirings, modulo\"{\i}ds, dio\"{\i}ds, pseudorings,
pseudomodules, band spaces over belts, idempotent mathematics). But first, one word of advise is in order.
Tropical conics are, of course, fairly well understood by
experts (in terms of combinatorics: secondary polytopes of
matrices, Gale dual spaces,  etc.), see \cite{Dickenstein}. Also, there exist
algorithms and computer programs to deal with them. Our
point is, nonetheless, that all of this can be done in
elementary terms, easily and fast, just by hand.

This paper  originated as an attempt to explain in full
detail and give proofs for all statements made in example
3.4 in \cite{Richter}.

Our polynomials will be either homogeneous in three
variables or non--necessary homogeneous in two variables.
To a degree--two tropical polynomial $p$, we associate a
point in the tropical plane and a triple of non--negative real numbers,
$s_{21}^+,s_{32}^+,s_{31}^+$, which completely determine  the
tropical conic  $\CC(p)$. These data are simply computed from $p$ and they
are all that is  needed to know in order to sketch $\CC(p)$. It is known
that the regular subdivision of the Newton polygon of $p$
determines the combinatorial type of $\CC(p)$ but, to our knowledge, nothing precise
has been said about the exact coordinates of the vertices
of $\CC(p)$.

There are two types (with several sub--types) of tropical conics: degenerate and non--degenerate ones. We explain
how they are classified according to the values of the
invariants $s_{21}^+,s_{32}^+,s_{31}^+$ and certain alternating sums $d_1,d_2,d_3$ of the $s_{ij}^+$'s. Degenerate (also called improper) tropical
conics are classified in theorem  \ref{thm:deg_con}.  It
turns out that pairs of tropical lines are degenerate
tropical conics, but the converse is not true. And non--degenerate (also called proper) tropical
conics are classified in theorem  \ref{thm:non_deg_con}, into \emph{one--point central} and \emph{two--point central} ones.

Given a degree--two tropical polynomial $p$, the values
$s_{21}^+,s_{32}^+,s_{31}^+$ can be arranged into a $3\times3$
symmetric non--negative real matrix denoted $\SHP(p)^+$. We characterize
tropical conics $\CC(p)$ having tropically singular
associated matrix $\SHP(p)^+$
 (corollary \ref{cor:singular}). These are
pairs of tropical lines and, surprisingly enough,
one--point central conics.

In the last section of the paper, we address the question
of irreducibility of degree--two  tropical polynomials, also in
elementary terms. We show that there exist irreducible
degree--two tropical polynomials giving rise to  pairs of
tropical lines.

 \m
Some of the results in this paper have already appeared in
\cite{Ansola}, while other are new.  %An earlier version of
%this paper is found in arXiv:math/0702143.
The idea of
considering shape matrices comes, somehow, from
\cite{Izhakian}. The values $s_{ij}$ come from \cite{Richter}.

Many  results in tropical algebra have been discovered since the late fifties
so that  the literature on this topic is vast. Some references are
the  books \cite{Baccelli, Cuninghame, Gunawardena, Zimmermann} and the  papers \cite{Akian, Butkovic, Cuninghame_B, Gaubert, Wagneur_F, Wagneur_M}.
The factorization problem for tropical polynomials in one
variable  has been
investigated in \cite{Kim}. The  tropical
version  of the  existence and uniqueness of a tropical
conic passing through five given points in the plane in
general position can be found in \cite{Richter}.

\m We would like to  thank the anonymous referee for pointing out a better way to present this material.
 %and to B.
%Bertrand, A. D{\'\i}az--Cano and
%P.D. Gonz{\'a}lez for some helpful discussions about the
%material presented here.

\section{Tropical conics}
\subsection{Tropical planes}

%Tropical conics are degree--two plane tropical curves and
%they can be defined in the affine setting or in the
%projective one.
\emph{Tropical geometry} arises when one works over $\T$,
the \emph{tropical semi--field}. By definition, $\T$ is the
set $\realamp$ endowed with two operations: $\max$ and $+$.
\emph{Tropical addition} is $\max$   and $+$ is taken as
\emph{tropical multiplication}. They are denoted $\oplus$
and  $\odot$, respectively. The neutral element for
tropical addition is $-\infty$  and zero is the neutral
element for tropical multiplication. It is noticeable that
$a\oplus a=a$, for $a\in\T$, that is, tropical addition is
\emph{idempotent}. However, there does not exist an inverse
element,
    with respect to $\oplus$, for  $a\in\T$. This
    is all that  $\T$ lacks in order to be  a field.

\m $\R_{\ge0}$ will denote the set of non--negative real
numbers. For $a\in\T$, we will set $a^+=\max\{a ,
0\}=a\oplus 0$, the \emph{non--negative part of} $a$.
For a matrix $A$,
the matrix obtained by replacing every entry $a$ of $A$ by
$a^+$ will be denoted $A^+$.  For a  polynomial $P$,
the polynomial obtained by replacing every coefficient $a$
of $P$ by $a^+$ will be denoted $P^+$.

\m The  \emph{tropical affine $2$--space} is  $\T^2$,
where addition and multiplication are defined
coordinatewise. It will be denoted $\T\A^2$. Here we can
define \emph{translations} %and \emph{homotheties}
in the standard way; every point $(t_1,t_2)\in\R^2$ defines the map: $(X,Y)\mapsto (X+t_1,Y+t_2)$.

\m  In the space
$\T^{3}\setminus\{(-\infty,-\infty,-\infty)\}$ we define
an equivalence relation $\sim $ by letting $(b_1,
b_2,b_{3})\sim (c_1, c_2,c_{3})$ if there exists

$\lambda\in\R$ such that
    $$(b_1+\lambda,  b_2+\lambda,b_{3}+\lambda)=
    (c_1, c_3,c_{3}).$$
The equivalence class of $(b_1, b_2,b_{3})$ is denoted
$[b_1, b_2,b_{3}]$. Its elements are obtained by adding
multiples of the vector $(1,1,1)$ to the point $(b_1,
b_2,b_{3})$. The \emph{tropical projective $2$--space},
$\T\P^2$, is the set of such equivalence classes. Notice
that, at least, one of the coordinates of any point in
$\T\P^2$ must be finite.

\m Points in $\T\A^2$ (resp. $\T\P^2$) having finite
coordinates will be called \emph{interior points}. The
rest of the points will be called \emph{boundary points}.
The \emph{boundary} of $\T\A^2$ (resp. $\T\P^2$) is the
union of its boundary points.  We will use $X,Y,Z$ as
variables in $\T\P^2$. Any permutation of the variables
$X,Y,Z$ provides a \emph{change of projective tropical
coordinates}. \emph{Translations} \label{dfn:translation} are also natural changes of
projective tropical  coordinates:  given $[t_1,t_2,t_3]\in \R^3$, the point
$[X,Y,Z]$ maps to  $[X',Y',Z']=[X+t_1,Y+t_2,Z+t_3]$. We may write
\begin{equation}
[X',Y',Z']=[X,Y,Z]\odot D, \qquad D=\left(\begin{array}{ccc}
t_1&-\infty&-\infty\\
-\infty&t_2&-\infty\\
-\infty&-\infty&t_3\\
\end{array}\right).
\end{equation}
A particular case is the \emph{tropical identity matrix}
$$I=\left(
\begin{array}{ccc}
0\\
-\infty&0\\
-\infty&-\infty&0\\
\end{array}\right).$$
Here, tropical matrix multiplication is defined in the usual way, but using $\oplus$ and $\odot$.

\m The plane $\T\P^2$ is covered by three copies of
$\T\A^2$ as follows. There exist injective maps
$$j_3:\T\A^2\to\T\P^2,\quad (x,y)\mapsto \ptop xy0,\qquad j_2:
\T\A^2\to\T\P^2,\quad (x,z)\mapsto \ptop x0z,$$
    $$j_1:\T\A^2\to\T\P^2,\qquad (y,z)\mapsto \ptop 0yz$$
     and
$\T\P^2=\im j_3\cup\im j_2\cup\im j_1.$ The complementary
set  of, say, $\im j_3$ is
$$\{\ptop
xy{-\infty}:x,y\in\T\}.$$ Moreover, we have
$j_3(x,x)=[x,x,0]=[0,0,-x]$, for $x\in\T$. This means that
 the coordinate axis $Z$  in $\T\P^2$ is transformed by $j_3^{-1}$ into
the usual line $X=Y$ in $\T\A^2$. The negative $Z$
half--axis in $\T\P^2$  corresponds to  the north--east
direction in $\T\A^2$. Similarly for $j_2,j_1$.

\m It is easy to check that the set of interior points
 of $\T\P^2$ equals the intersection
$\im j_3\cap \im j_2\cap \im j_1.$

\m For simplicity, we will  consider the Euclidean metric
in $\T\A^2$. Notice that the composite maps $j_l^{-1}\circ
j_k$ are NOT isometries, for $k, l\in\{1,2,3\}$, $k\neq
l$.

\m The projective tropical coordinates of a point in
$\T\P^2$ are not unique. In order to avoid this
inconvenience,  we choose a \emph{normalization}, that is
we fix a rule that allows us to have unique coordinates
for all (but perhaps a small set of)  points in $\T\P^2$,
according to this rule. For instance, setting the last
coordinate equal to zero is a normalization. We call it
the \emph{$Z=0$ normalization} and  say that \emph{we work
in} $Z=0$. %For any given normalization, there exist a few
%points in $\T\P^2$ which do not admit it. This is the case
%for $\ptop ab{-\infty}$, for any $a, b\in\T$ and the $Z=0$
%normalization.
To consider the $Z=0$ normalization is the same thing as
passing to the affine tropical plane, via $j_3$. Other
possible normalizations are setting $Y=0$, or setting
$X=0$, or setting $X,Y,Z$ all non--negative and, at least,
one equal to zero, or setting $X+Y+Z=0$, etc.

\subsection{Tropical conics  are trees}\label{subsec:con_are_trees}
A tropical polynomial is a  tropical sum of tropical monomials.
For instance, a tropical homogeneous degree--two polynomial in the variables $X,Y,Z$ is
$P(X,Y,Z)$ $$=\poltrohdos a$$
$$=\maxxi a.$$ For us, degree--two means that the Newton polygon of $P$ is the triangle
determined by the points $(2,0),(0,2),(0,0)$; in other words, that
 $\sucespar a \in\T$ but $a_{11},a_{22},a_{33}\in\R$.
The \emph{tropical conic} $\CC(P)$ given by $P$ is, by
definition, the set of points in $\T\P^2$ where the
\emph{maximum is attained, at least, twice.}
 A simple computer program (done in MAPLE, for instance) may be
used in order to sketch this conic, say in $Z=0$. But we want to show
that one can easily sketch
$\CC(P)$ without a computer! Indeed, it is well--known that $\CC(P)$ is a tree,
see \cite{Gathmann, Mikhalkin_E, Mikhalkin_T, Richter} and so, all we need to compute is the coordinates of its vertices.

So let us
recall here some facts about \emph{trees}; see \cite{Foulds,
Harary} for details. A \emph{graph} $G$ is an ordered pair
$(V,E)$, where $V$  is a finite set of points, called
\emph{vertices}  of $G$, and $E$  is a set of
cardinality--two subsets of $V$. The elements of $E$ are
called \emph{edges} of $G$. The edge joining vertices $u,w$
will be denoted $\overline{uw}$. The  \emph{degree}
 of a vertex $w$ of $G$ is the number of edges of $G$
incident with $w$.  %%%; it will be denoted $\deg(w)$.
Degree--one vertices are called \emph{pendant vertices} and
edges incident to pendant
vertices are called \emph{pendant edges}. %%Let us write
%%$\deg(G)= \max_{v\in V}\deg(v)$; this is the \emph{degree
%%of the graph} $G$. NO LO USO

A \emph{tree} is a connected graph without cycles. %%In a
%%tree the number of edges is one less than the number of
%%vertices.
A tree  $G=(V,E)$ naturally carries a \emph{discrete
metric}; it is the function
    $\dd:V\times V\to \N$,
where $\dd(u,w)$ is the least number of edges to be passed
through when going from $u$ to $w$.  If $\dd(u,w)=1$, we
say that $u,w$ are \emph{consecutive vertices}. The
\emph{eccentricity} of a vertex $w$ is
    $\e(w)=\max_{u\in V}\dd(u,w)$ and the \emph{radius}
of the graph $G$ is
    $\rr(G)=\min_{w\in V}\e(w).$
A vertex $w$ in  $G$ is   \emph{central} in $G$ if
$\e(w)=\rr(G)$ and the \emph{center} of $G$ is the set of
all  central points in $G$. It is known that \emph{every
tree has a center and it consists either of just one vertex
or two consecutive vertices.} This explains the names
\emph{one--point central} and \emph{two--point central conics}, given below in theorem \ref{thm:non_deg_con}.
%%p.~\pageref{name} and \pageref{namee}. \label{name2}

A tropical projective plane curve $\CC$ of
degree $d>0$ is a weighted tree of a very particular sort.
Each vertex of $\CC$ is determined by its tropical
projective coordinates. The pendant vertices of $\CC$ are
precisely the points in $\CC$ which lie on the boundary of
$\T\P^2$. There are $3d$ such vertices, counted with
multiplicity. They are grouped in 3 families of $d$
vertices each: $d$ vertices have the $X$ (resp. $Y$) (resp.
$Z$) coordinate equal to $-\infty$. Every pendant edge in
$\CC$ has infinite length. There are $3d$ such edges,
counted with multiplicity,  and they  are grouped in
3 families of $d$ edges each. %%, one such family for each
%%coordinate direction.
The rest of the edges in $\CC$ have finite lengths. Edges
in $\CC$ may carry a multiplicity, which is a natural
number, no greater than $d$. The multiplicity of a vertex
is deduced from the multiplicities of the edges incident to
it.

 \m A tropical projective  plane  curve $\CC$ can
be represented in $Z=0$ (or in $Y=0$ or $X=0$). More
precisely, this means that we  represent $j_3^{-1}(\CC)$
(and still denote it $\CC$) (or $j_2^{-1}(\CC)$ or
$j_1^{-1}(\CC)$) in $\T\A^2$. Say, we choose to work in
$Z=0$. Then the \emph{slope} \label{rem:slope} of every
edge of finite length in $\CC$ is a rational number and at
each non--pendant vertex $w$ the \emph{balance condition}
holds. This means that $\sum_{j=1}^s\lambda_je_j=0$, where
$\suces u1s$ are all the vertices in $\CC$ consecutive to
$w$, $\suces\lambda 1s\in\N$ are the weights of the edges
$\overline{wu_1},\ldots \overline{wu_s}$ and $\suces
e1s\in\Z^2$ are the primitive integral vectors at the point
$w$ in the directions of such edges.

\subsection{Matrices and points associated to a tropical degree--two polynomial}

Let
$$P=\maxxi a$$ be a homogeneous tropical
polynomial of degree two. As explained in
subsection \ref{subsec:con_are_trees}, the tropical conic $\CC(P)$ has six pendant
edges, counted with multiplicities. These multiplicities
are either one or two. Without loss of generality, we may
work in $Z=0$. Then $\CC(P)$ has two pendant edges in the
west direction, two in the south direction and two in the
north-east direction, all counted with multiplicity.
In order to sketch the conic $\CC(P)$ we must determine the
non--pendant vertices of $\CC(P)$. We will see that there are  four such points, at most.

Just like in usual geometry, to $P$ we associate the symmetric matrix $A(P)=(a_{ij})$,
bearing in mind that we need not divide the coefficients of mixed terms by two, since tropical addition is idempotent.
Conversely, to such a matrix $A$, we can associate a polynomial $P(A)$ and, eventually, a tropical conic $\CC(A)$.

Most matrices considered in this paper are $3\times 3$ and have entries in $\T=\realamp$
but their diagonal entries belong to $\R$ (the only exception appears in the definition of tropical determinant) and  are
symmetric. Therefore, we only
write their lower triangular parts.
%We will arrange
%subscripts of the entries in an unusual way, taken from
%\cite{Richter}, and all subscripts will be taken modulo 6.
To the symmetric matrix  $$A=\mat a$$ we associate the diagonal matrix
$$D=D(A)=\diag a,$$ which corresponds to a translation of coordinates, as we have seen in p. \pageref{dfn:translation}.
The
tropical inverse matrix of $D$  is
obtained by negating the signs of its diagonal entries. Obviously, it corresponds to the inverse translation.
%Then, we consider the change of coordinates (a translation)  given by
%$[X,Y,Z]\mapsto  [X,Y,Z]\odot D(A)=[X+\dfrac{a_{11}}{2},Y+\dfrac{a_{22}}{2},Z+\dfrac{a_{33}}{2}]$
We  define the \emph{shape matrix} associated to $A$ as $S=\SHP(A)=D^{\odot -1}\odot A \odot D^{\odot -1}$.
Clearly, \textbf{the shape matrix corresponds to the given conic $\CC(P)$, after translation}.
It is crucial and easy to check that \emph{the shape matrix  $S=(s_{ij})$ is symmetric and has zero diagonal entries.}
The remaining entries of $S$  are related to $A$ by the following formulas:
\begin{equation}\label{eqn:s_en_a}
2s_{21}=2a_{21}-a_{11}-a_{22},\quad 2s_{32}=2a_{32}-a_{22}-a_{33},\quad
2s_{31}=2a_{31}-a_{33}-a_{11}.
\end{equation}
Therefore the shape matrix  is
 \emph{invariant}, in the sense that  it does not change if $A$
is replaced by $A=\alpha+U$, for any $\alpha\in\R$, where
$U$ denotes the $3\times 3$ matrix all whose entries are
one.
Notice  also that  the matrices $A$ and $S$ are/are not simultaneously  real.
 Back to the polynomial $P$, let  $\SHP(P)$,  $D(P)$  denote the polynomials  associated to the  matrices $S$ and $D$.

\m  The \emph{tropical determinant} of an arbitrary $3\times 3$
 matrix $A=(a_{ij})$
 is defined as $$|A|_{trop}=\max_{\sigma\in S_3}\{a_{1\sigma(1)}+a_{2\sigma(2)}+a_{3\sigma(3)}\},$$
where $S_3$ denotes the permutation group in $3$ symbols. A
matrix is  \emph{tropically singular} if the \emph{maximum
in the tropical determinant is attained, at least, twice}.
%For a matrix as in (\ref{eqn:matrix}),
For the matrices above we have
$$2|D|_{trop}=a_{11}+a_{22}+a_{33}$$
and $D$ is tropically non--singular. Moreover, $A$  and $S$  are/are not simultaneously  tropically singular, because
$$\sum_{i=1}^{3}a_{i\sigma(i)}=\sum_{i=1}^{3}s_{i\sigma(i)}+2|D|_{trop},$$ for all $\sigma\in S_3$. The tropical
determinants of $S$  and $A$  are equal to
$\max\{0,s,s,2s_{21},2s_{32},2s_{31}\}$ and
$2|D|_{trop}+|S|_{trop}$, respectively,
where
\begin{equation}\label{eqn:s}s=s_{21}+s_{32}+s_{31}.
%%=2(-a_{11}+a_{21}-a_{22}+a_{32}-a_{33}+a_{31}),
\end{equation}

\begin{lem}
$\SHP(\SHP(A))=\SHP(A).$
\end{lem}

\begin{proof}
It follows from the formulas (\ref{eqn:s_en_a}).
\end{proof}

In the following, we  assume
$A=\SHP(A)$ (or equivalently, $P=\SHP(P)$), meaning  that $a_{11}=a_{22}=a_{33}=0$ and $a_{ij}=s_{ij}$, if $i\neq j$.
Now, the next crucial lemma tells us that  \emph{the matrices $A$
 and $A^+$ give rise  to the same tropical conic.}

\begin{lem}\label{lem:P_P+_equal} If
$P=\SHP(P)$, then $\CC(P)=\CC(P^+)$.
\end{lem}

\begin{proof} By hypothesis,
$P=\max\left\{
2X, 2Y, 2Z, {s_{21}}+X+Y , {s_{32}}+Y+Z ,
 {s_{31}}+X+Z
 \right\}.$
If
$-\infty\le s_{21}\le0$ then $P^+=\max\left\{
2X, 2Y, 2Z, X+Y, {s_{32}^+}+Y+Z ,
{s_{31}^+}+X+Z
 \right\}.$
It is obvious that
$$ \max\{2X ,2Y, {s_{21}}+X+Y\}=\max\{2X ,2Y\}=\max\{2X ,2Y,X+Y\}.$$
Moreover, these three maxima are attained at least twice at exactly the same points in $\R^2$. Therefore,  the term $X+Y+{s_{21}}$
is irrelevant in $P$, as far as $\CC(P)$ is concerned.

We can reason similarly with $s_{32},s_{31}$,
and thus conclude that the tropical conics $\CC(P),\CC(P^+)$ are equal.
\end{proof}

In the former paragraphs, we have reduced the study of tropical conics to the case $A=\SHP(A)^+$, a non--negative real matrix.
Now,   what does  such a tropical conic $\CC(A)$ look like, say in $Z=0$?
To answer this question, we define  the points $v^1(A),v^2(A),v^3(A)$ which arise from
 the rows of  $A$:
$$v^1(A)=\hh s,\ v^2(A)=\vv s,\ v^3(A)=\ii s$$
and one more point $v^0(A)$ by
$$v^0(A)=[s_{32},s_{31},s_{21}].$$ The points will be denoted $v^0,v^1,v^2,v^3$, for short.

\begin{lem}\label{lem:triedro} Suppose  $A=\SHP(A)^+$. Then,  in $Z=0$,
\begin{enumerate}
\item the segment $\overline{v^0v^1}$ is parallel to the $X$  axis,
\item the segment $\overline{v^0v^2}$ is parallel to the $Y$  axis,
\item the segment $\overline{v^0v^3}$ is parallel to the line $X=Y$.
\end{enumerate}
\end{lem}

\begin{proof} Taking differences, we have $v^1-v^0=[-s_{32}, -s_{21}-s_{31}, -s_{21}-s_{31}]$ and the coordinates of this point  in
$Z=0$ are $(-s_{32}+s_{21}+s_{31},0)$. The rest is similar: $v^2-v^0=(0, -s_{31}+s_{32}+s_{21})$ and $v^3-v^0=(-s_{31}-s_{32}+s_{21}, -s_{31}-s_{32}+s_{21})$.
\end{proof}

Notice how the lengths  of the segments $\overline{v^0v^j}$ depend on alternating sums of the entries of the matrix $A=\SHP(A)^+$.
More precisely, set
\begin{equation}
\left(\begin{array}{c} d_1\\d_2\\d_3\\
\end{array}\right)=\left(\begin{array}{ccc} 1&-1&1\\ 1&1&-1\\ -1&1&1\\ \end{array}\right) \left(\begin{array}{c} s_{21}\\ s_{32}\\ s_{31}\\ \end{array}\right) \label{eqn:d_en_s}
\end{equation}
(in terms of the ordinary matrix multiplication). Hence
\begin{equation}\label{eqn:s_en_d}
s_{ij}=\dfrac{d_i+d_j}{2},\quad  i\neq j.
\end{equation}
The length of  $\overline{v^0v^j}$ is $|d_j|$, for  $j=1,2$,  and the  length of  $\overline{v^0v^3}$ is $\sqrt2|d_3|$
(the factor $\sqrt5$ is due to our choice of normalization $Z=0$).
Moreover, the angle $\angle v^1v^0v^2$ is $\dfrac{\pi}{2}$. %; the signs of $d_1$ and $d_2$ determine
% the triangle $v^1,v^0,v^2$ further.
In addition,  $\angle v^1v^0v^3$ is
 $\frac{3\pi}{4}$ (resp. $\frac{\pi}{4}$)  if $d_1d_3>0$
 (resp. $d_1d_3<0$). Notice that the vertices $v^1,v^0,v^3$ determine a  right triangle
 in $Y=0$. Similarly,  the vertices $v^2,v^0,v^3$ determine a  right triangle
 in $X=0$.
%Also the vertices $v^2, v^3,
% v^0$ determine
% an usual triangle in $Z=0$, having at $v^0$ an angle of
% $\frac{3\pi}{4}$ (resp. $\frac{\pi}{4}$)  if $d_2d_3>0$
% (resp. $d_2d_3<0$). These triangles become right triangles
% in $X=0$.

\begin{lem}\label{lem:positivo2} If $A=\SHP(A)^+$, then
$d_j$ is negative for, at most, one $j\inimp$. % And if
%$\sucespar f$ are all $\ge 0$, then  $d_j$ is $\le0$, for,
%at most, one $j\inimp$. \qed
\end{lem}

\begin{proof}  Suppose $d_1<0$. By the hypothesis and relations (\ref{eqn:s_en_d}), $0\le d_1+d_2$
and $0\le d_1+d_3$, whence $0<-d_1\le d_2$ and $0<-d_1\le
d_3$.  The other cases are similar.
\end{proof}

\begin{lem}\label{lem:singular} If $A=\SHP(A)^+$, then
the following are equivalent:
\begin{enumerate}
\item $A$ is tropically singular,
\item the maximum of $\sucespar s$ is no greater than the
sum of the other two,
\item $\sucesimp d$ are all non--negative.
\end{enumerate}
\end{lem}

\begin{proof}
Equivalence between 2. and 3. follows from (\ref{eqn:d_en_s}).  We show that 1. and 2. are equivalent.
Note that only one of the numbers $2s_{21},2s_{32},2s_{31}$ can be greater than $s=s_{21}+s_{32}+s_{31}$,
and therefore, the maximum in $\max\{0,s,s,2s_{21},2s_{32},2s_{31} \}$ is
attained twice if and only if it is equal to $s$.  Now note that this happens if and only if 2. is satisfied.
%Without loss of generality, we may suppose that
%$s_{31}=\max\{s_{21},s_{32},s_{31}\}$. Then $s_{31}\le s_{21}+s_{32}$ if and
%only if $2s_{31}\le s$ if and only if $s=\max\{0,s,s,
%2s_{21},2s_{32},2s_{31}\}$   if and only if Equivalence  with $d_1\ge0$,
%$d_2\ge0$, $d_3\ge0$ follows from (\ref{eqn:d_en_s}).
\end{proof}

Any tropical conic  $\CC$  has some non--pendant vertices.
These are the   points in $\CC$ where the maximum is attained, at least, three times.
 %A tropical conic   is \emph{non--degenerate} if it has the maximal number of non--pendant vertices. %This number is four.
If  $\CC$  has more than two non--pendant vertices,  let us consider  two consecutive ones $u^1,u^2$. If these points come together,   a new tropical conic $\CC'$ arises. Obviously,  if $\CC$
has parallel  pendant edges $e^1,e^2$  such that $e^j$ is
incident to $u^j$, then   $e^1$ is a
pendant edge with multiplicity two in $\CC'$.
\label{dfn:degeneration}
Let $\CC'$ be a tropical conic which can be obtained from
$\CC$ by successively collapsing one or more pairs of
consecutive  non--pendant vertices. Then we will say that
\emph{$\CC'$ is a degeneration of $\CC$}. Such a conic $\CC'$ is called \emph{degenerate}.
%\end{dfn}

\m Now we get our two main theorems. In page \pageref{rem:comentario} we explain why theorem
\ref{thm:non_deg_con} deals with \emph{non--degenerate tropical conics} while theorem \ref{thm:deg_con} classifies  \emph{degenerate tropical conics}.

In the second part of the  following theorem, superscripts work modulo 3, and $t^{i,j}$ stands for the point in $\T\P^2$ whose  $i$--th coordinate is  $-2s_{ij}$ and the rest are null.

\begin{thm}\label{thm:non_deg_con}
Let $A=\SHP(A)^+=(s_{ij})$. Suppose that  $s_{ij}>0$ for all $i\neq j$ and $d_j\neq0$, for $j=1,2,3$. Then the following mutually exclusive cases arise, for the tropical conic $\CC=\CC(A)$.
\begin{enumerate}
\item \emph{One--point central conic}. If $d_1, d_2, d_3$ are all positive, then  $\CC$ has four non--pendant vertices; these are $v^1,v^2,v^3$ and $v^0$.
\item \emph{Two--point central conic}.
If $d_j<0$ for some $j\in\{1,2,3\}$, then  $\CC$ has four non--pendant vertices; these are $v^{j-1},v^{j+1}, w^{j-1}=v^{j-1}+t^{j-1,j}$ and $w^{j+1}=v^{j+1}+t^{j+1,j}$.
%If $d_3<0$, then  $\CC$ has four non--pendant vertices; these are $v^1,v^2, w^1=v^1+[-2s_{31},0,0]$ and $w^2=v^2+[0, -2s_{32},0]$. Similarly if $d_1<0$ or $d_2<0$.
\end{enumerate}
\end{thm}

\begin{figure}[H]
\centering
\includegraphics[width=13.5cm, keepaspectratio]{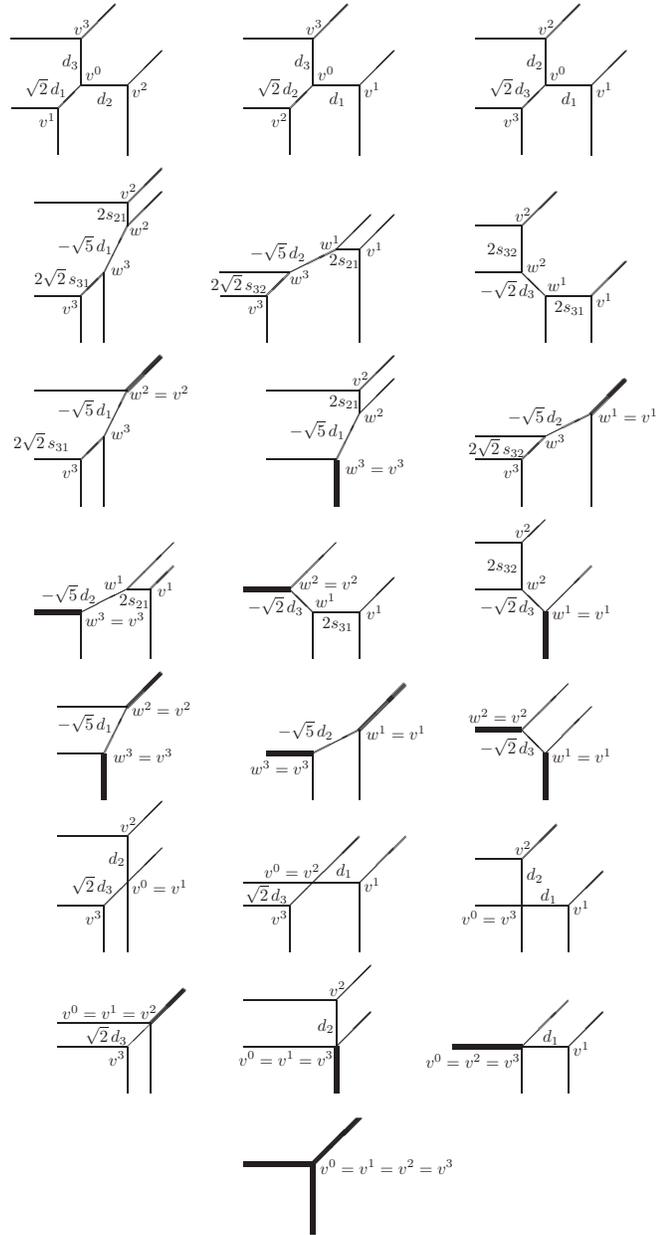}\\
\caption{Tropical conics. Line 1 is occupied by one--point central conics, line 2 is occupied by two--point central conics, lines 3 to 8 are occupied by degenerate conics, where lines 6 to 8 are occupied by pairs of  lines.} \label{fig:conicas}
\end{figure}

\begin{proof} We may assume that $d_1>0$ and $d_2>0$ by  a permutation of variables and lemma \ref{lem:positivo2}.
For simplicity, let us work in $Z=0$ and let us evaluate $$P=\max\{2X,2Y,0,s_{21}+X+Y,s_{31}+X,s_{32}+Y\}$$ in $v^1=(s_{31}, s_{31}-s_{21})$ and $v^2=(s_{32}-s_{21}, s_{32})$.
Using that $d_1>0$ and $d_2>0$, we obtain that
$$\max\{2s_{31}, 2(s_{31}-s_{21}),0,2s_{31},2s_{31}, d_3\}=2s_{31}$$
$$\max\{2(s_{32}-s_{21}), 2s_{32},0,2s_{32}, d_3,2s_{32}\}=2s_{32}$$
both attained three times. This means that $v^1$ and $v^2$ are non--pendant vertices of $\CC$.
Now we evaluate $P$ in $v^3=(-s_{31},-s_{32})$ and $v^0=(s_{32}-s_{21},s_{31}-s_{21})$ and obtain
$$\max\{-2s_{31}, -2s_{32},0,-d_3, 0,0\}$$
$$\max\{2(s_{32}-s_{21}), 2(s_{31}-s_{21}),0,d_3,d_3,d_3\}.$$
It follows that $v^3$ and $v^0$ are also non--pendant vertices of $\CC$, if $d_3>0$ and, no further non--pendant vertices of $\CC$ arise, by symmetry in the variables; see figure \ref{fig:interior}, right. The  center of $\CC$ is $v^0$ and we say that $\CC$ is a \emph{one--point central conic}.
Six
pendant edges hang from the $v^1,v^2,v^3$ as explained in  subsection \ref{subsec:con_are_trees}, completing the picture of $\CC$;  see
figure \ref{fig:conicas} line 1, column 3.
If we work in $X=0$, (resp.
$Y=0$) we obtain other representations of $\CC$; see
figure \ref{fig:conicas} line 1, column 1 (resp. column 2).
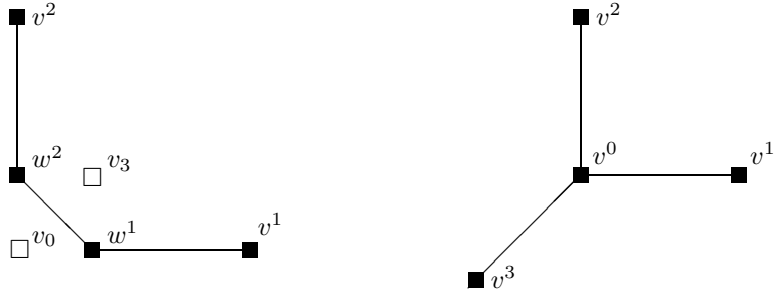
\begin{figure}[H]
\setlength{\unitlength}{1cm} %Unidad para gráficas.
%\pagestyle{empty}
%\setlength{\topmargin}{-3cm}
%\begin{document}
 \begin{center}

  \begin{picture}(12,5)(0,0)
%Cónica 1-punto-centrada
    \put(7,0){\line(1,1){1.5}}
    \put(8.5,1.5){\line(0,1){2}}
    \put(8.5,1.5){\line(1,0){2}}
    %puntos
    \put(7,0){\rule{2mm}{2mm}}
    \put(8.4,3.5){\rule{2mm}{2mm}}
    \put(8.4,1.4){\rule{2mm}{2mm}}
    \put(10.5,1.4){\rule{2mm}{2mm}}
    %Etiquetas
    \put(7.3,0){$v^3$}
    \put(8.65,1.65){$v^0$}
    \put(10.75,1.65){$v^1$}
    \put(8.7,3.5){$v^2$}
%Cónica 2-puntos-centradas
    \put(1,1.5){\line(0,1){2}}
    \put(2,0.5){\line(1,0){2}}
    %Pendiente -1
    \put(1,1.5){\line(1,-1){1}}
    %puntos
    \put(0.9,1.4){\rule{2mm}{2mm}}
    \put(0.9,3.5){\rule{2mm}{2mm}}
    \put(1.9,0.4){\rule{2mm}{2mm}}
    \put(4,0.4){\rule{2mm}{2mm}}
    %Etiquetas
    \put(4.2,0.7){$v^1$}
    \put(1.2,3.5){$v^2$}
    \put(2.2,0.6){$w^1$}
    \put(1.2,1.6){$w^2$}
    %vo y v3
    \put(0.9,0.4){$\square$}
    \put(1.87,1.35){$\square$}
    \put(1.2,0.6){$v_0$}
    \put(2.2,1.6){$v_3$}
    \end{picture}
   \end{center}
\caption{Non--pendant vertices: cases $d_3<0$ and $d_3>0$.} \label{fig:interior}
\end{figure}
Now, if   $d_3<0$,  we consider $w^1=v^1+[-2s_{31},0,0]$ and $w^2=v^2+[0, -2s_{32},0]$. Working in $Z=0$ and evaluating $P$ in $w^1=(-s_{31}, s_{31}-s_{21})$ and $w^2=(s_{32}-s_{21}, -s_{32})$ we get
$$\max\{-2s_{31}, 2(s_{31}-s_{21}),0,0,0, d_3\}=0$$
$$\max\{2(s_{32}-s_{21}),-2s_{32}, ,0,0, d_3,0\}=0$$
both attained three times.
It follows that $w^1$ and $w^2$ are non--pendant vertices of $\CC$ (in addition to $v^1$ and $v^2$), if $d_3<0$. No more non--pendant vertices of $\CC$ arise also in this case.
In particular,
 $v^3, v^0$ are NOT vertices in $\CC$, if $d_3<0$; see figure \ref{fig:interior}, left. The  center of $\CC$ consists of $w^1$ and $w^2$ and we say that $\CC$ is a \emph{two--point central conic}. Six pendant edges of
$\CC$  hang from $v^1,v^2,w^1,w^2$. Such a tropical conic is represented
in figure \ref{fig:conicas} line 2, column 3.

If $d_1<0$ or $d_2<0$, other two--point central conics are obtained, and they are represented in figure \ref{fig:conicas} line 2, columns 1 and 2.
Notice  that a factor  $\sqrt5$
appears in the length of edges  of slope
$\frac{1}{2}$ or 2, due to our choice of Euclidean metric.
\end{proof}

\begin{cor}\label{cor:distancia_aristas_no_aco}
Let $A=\SHP(A)^+=(s_{ij})$. Suppose that  $s_{ij}>0$ for all $i\neq j$, $d_1>0$, $d_2>0$ and $d_3\neq0$.
Then in $Z=0$, the tropical conic $\CC=\CC(A)$  has two
different pendant edges in the north--east direction (resp.
 west direction) (resp. south direction) and they are
$\sqrt{d_1^2+d_2^2}$ (resp. $2s_{32}$) (resp. $2s_{31}$) apart.
\end{cor}
\begin{proof} The previous theorem applies and the statement follows from the equalities
(\ref{eqn:s_en_d}).
\end{proof}

Notice that $\sqrt{d_1^2+d_2^2}$ tends to zero if and only
if $2s_{21}=d_1+d_2$ tends to zero.

\m Suppose that $A=\SHP(A)^+$ and $s_{ij}>0$, for $i\neq j$ and $d_j\neq0$, for $j=1,2,3$. Then the degenerations of the tropical conic $\CC(A)$ arise by letting  $s_{ij}=0$  or $d_j=0$ for some indices. We have the following classification theorem.

\begin{thm}%[Degenerated tropical conics]
\label{thm:deg_con}
If  $A=\SHP(A)^+$ and $s_{ij}=0$  or $d_j=0$ for some indices, then the following  cases arise (up to a permutation of variables) for the tropical conic $\CC(A)$:
\begin{enumerate}
\item  $s_{21}>0, s_{32}>0, s_{31}=0, d_1>0, d_2>0$ and $d_3<0$.
\item  $s_{21}>0, s_{32}=s_{31}=0, d_1>0, d_2>0$ and $d_3<0$.
\item \emph{Double tropical line}. $s_{21}=s_{32}=s_{31}=0$ (equivalently, $d_1=d_2=d_3=0$ or, yet equivalently, $d_3=s_{32}=s_{31}=0$).
\item \emph{Pair of tropical lines}. $s_{21}>0, s_{32}>0, s_{31}>0, d_1>0, d_2>0$ and $d_3=0$.
\item \emph{Pair of tropical lines}. $s_{21}>0, s_{32}=0, s_{31}>0, d_1>0$ and $d_2=d_3=0$.
\end{enumerate}
\end{thm}

\begin{proof}
 \begin{enumerate}
 \item This situation arises when $v^1$ and $w^1$ collapse, in a two--point central conic.
 \item This situation arises when, in addition to the former,  $v^2$ and $w^2$ collapse, in a two--point central conic.
 \item This situation arises when $v^1,w^1, v^2$ and  $w^2$  all collapse to one point, in a two--point central conic.
 It also arises when $v^j$ all collapse to one point, for $j=0,1,2,3$, in a one--point central conic.
 \item This situation arises when $w^1$ and $w^2$ collapse, in a two--point central conic. It also arises when  $v^0$  and $v^3$ collapse, in a   one--point central conic.
 \item This situation arises when  $v^2, w^2$ and $w^1$ all collapse, in a two--point central conic.  It also arises when  $v^0, v^2$  and $v^3$ collapse, in a   one--point central conic.
 \end{enumerate}
 \end{proof}

These conics are represented in figure
\ref{fig:conicas}, lines 3 to 8, where a thick segment
represents a multiplicity--two edge.

\m Let us summarize.  Up to translation, tropical conics are determined by a non--negative real matrix $S^+=(s_{ij})$ with diagonal zero entries. We have gone through all the possibilities for the $s_{ij}$, in the two theorems above. This  means that no more tropical conics do exist!  Therefore,  theorem \ref{thm:non_deg_con} classifies \emph{non--degenerate tropical conics}, while  theorem \ref{thm:deg_con}  classifies \emph{degenerate tropical conics}.\label{rem:comentario}

\m \label{proc:draw} A \textbf{procedure to
sketch, say in $Z=0$, the tropical conic $\CC(P)$ defined by
an arbitrary homogeneous degree--two polynomial} $P$ is the following
\begin{itemize}
\item From $P$,  compute the matrices $A$ and
$S^+=\SHP(A)^+$.% and the translation point
%$\tau=\dfrac{1}{2}\left(a_{33}-a_{11},a_{33}-a_{22} \right)\in\R^2$.

\item Sketch the conic $\CC(S^+)$, according to the classification given by the theorems above and
translate this conic to the point $\dfrac{1}{2}\left(a_{33}-a_{11},a_{33}-a_{22} \right)$ in $\R^2$ to obtain $\CC(P)$.
\end{itemize}

The following are all direct consequences of our
discussion.

\begin{cor}
A tropical conic is non--degenerate if and only if it is not the union of
two  tropical lines and all of its pendant edges have
multiplicity one.\qed
\end{cor}

\begin{cor} [Pairs of tropical lines] \label{cor:pair}
For a tropical conic $\CC=\CC(A)$, the following statements are equivalent:
\begin{itemize}
\item $\CC$ is a pair of
 lines,
 \item  $\sucesimp d$ are all non--negative and, at
 least,  one $d_j$ equals zero,
 \item the maximum of $\sucespar{s^+}$ equals the
sum of the other two,
\item   $v^0\in\{v^1,v^2,v^3\}$ for the matrix $\SHP(A)^+$.\qed
\end{itemize}
 \end{cor}

Notice that the number of different pendant edges in a
pair of tropical lines is six, five or three.
Pairs of tropical lines are
represented in figure \ref{fig:conicas}, lines 6 to 8.

\begin{cor}\label{cor:singular}
A tropical conic $\CC=\CC(P)$ has tropically singular
associated matrix $\SHP(P)^+$
if and only if $\CC$ is either a
pair of tropical lines or  a
one--point central conic.
\end{cor}

\begin{proof}
This follows from lemma \ref{lem:singular}, part 1 of theorem \ref{thm:non_deg_con}   and corollary \ref{cor:pair}.
\end{proof}

\begin{figure}[H]
\centering
\includegraphics[width=8cm, keepaspectratio]{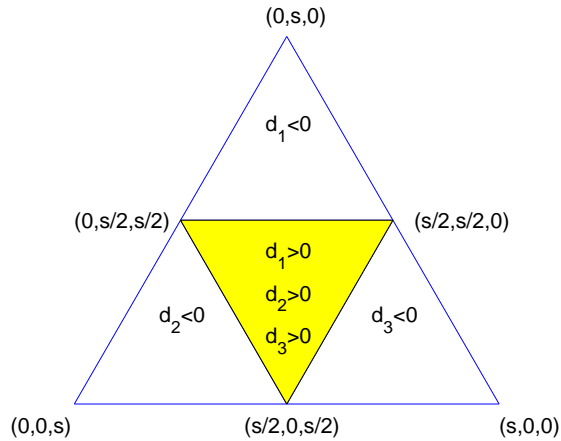}\\
\caption{Section of octant $\OO$ given by
$s_{21}^++s_{32}^++s_{31}^+=s$.} \label{fig:section}
\end{figure}

A tropical conic $\CC(P)$ is determined by a triple $(s_{21}^+,s_{32}^+,s_{31}^+)$ of real non--negative numbers and
any row of  the matrix $A=A(P)$. The null triple corresponds to a double tropical line.
 Let $(s_{21}^+,s_{32}^+,s_{31}^+)\neq(0,0,0)$ be the
coordinates  of a point in the non--negative octant
$\OO=\R_{\ge0}^3$. In figure
\ref{fig:section} we see the plane section of $\OO$ given
by $s_{21}^++s_{32}^++s_{31}^+=s$, for some positive $s$. According to corollary \ref{cor:singular}, tropical
conics having tropically singular matrix $\SHP(A)^+$
correspond to the shaded  closed triangle, the boundary of
which corresponds to pairs of lines. Other degenerate tropical
conics correspond to the
boundary of the section.

\m It is known that every  balanced weighted tree is a
tropical curve, see
\cite{Gathmann, Mikhalkin_D, Mikhalkin_E}. When $d=2$, here
is a \label{proc:poly} \textbf{procedure to find a
defining polynomial $P$ for  a balanced graph $\CC$}.
\begin{itemize}
\item From the edges of $\CC$, compute the values $\sucespar {s^+}$ %, \sucesimp
%d$ and the matrix $\SHP(A)^+$. C
and classify $\CC$ (degenerate or non--degenerate and type).

\item From the vertices of $\CC$, compute  symmetric matrices $A$ and $\SHP(A)$, using as many unknowns as
necessary.

\item Solve for the unknowns, according to the classification.
\end{itemize}
\m
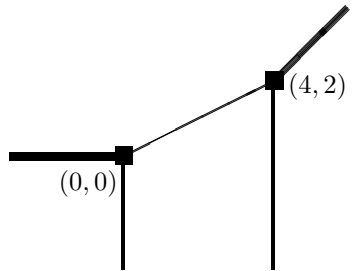
\begin{figure}[H]
\setlength{\unitlength}{5mm} %Unidad para gráficas.
\setlength{\topmargin}{-2cm}

 \begin{center}
  \begin{picture}(10,6)(0,0)
    \linethickness{3pt}
    \put(0,3){\line(1,0){3.04}}
    \linethickness{1pt}
    \put(3,0){\line(0,1){3}}
    \put(7,0){\line(0,1){5}}
    %puntos
    \put(2.8,2.8){\rule{2.3mm}{2.3mm}}
    \put(6.8,4.8){\rule{2.3mm}{2.3mm}}
    %pendiente 1/2
    \put(3,3.02){\line(2,1){4}}
    \put(3.02,3){\line(2,1){4}}
    %ramas infinitas
    \multiput(6.92,5.02)(0.02,-0.02){6}{\line(1,1){2}}
    %Etiquetas
    \put(1.3,2.1){$(0,0)$}
    \put(7.4,4.7){$(4,2)$}

  \end{picture}
 \end{center}
 \caption{A weighted tree $\CC$ in $Z=0$.} \label{fig:conica02}
\end{figure}

\begin{ex}\label{ex:grafo_equilibrado}
In  $Z=0$, let the weighted tree $\CC$ in figure
\ref{fig:conica02}
be given. Here thick  segments
represent edges of  multiplicity two.  The non--pendant vertices of $\CC$ are
$v^1=(4,2)=[4,2,0], v^3=(0,0)=[0,0,0]$ and
the balance condition is satisfied at both. Indeed, at $v^1$ the primitive vectors are
$(1,1),(0,1),(-2,-1)$ and $2(1,1)+(0,1)+(-2,-1)=(0,0)$. Similarly, for $v^3$. Therefore,
this tree corresponds to a tropical conic.
It is a degenerate tropical conic (not a pair of lines) and, by corollary \ref{cor:distancia_aristas_no_aco} and theorem \ref{thm:deg_con},
$s_{21}^+=s_{32}^+=0$ and $s_{31}^+=2$.
Then
$s_{21},s_{32}$ are non--positive and $s_{31}=2$. We  fill the negated coordinates of $v^1,v^3$
into the rows of a symmetric matrix $A$ and compute  the matrices $D=D(A)$ and $\SHP(A)=D^{\odot -1}\odot A \odot D^{\odot -1}$ obtaining:
$$A=\left(
\begin{array}{ccc}
-4\\
-2&a_{22}\\
0&0&0\\
\end{array}\right), \qquad
\SHP(A)=\left(
\begin{array}{ccc}
0\\
-a_{22}/2&0\\
2&-a_{22}/2&0\\
\end{array}\right),$$
 for some $a_{22}\in\R$. Therefore  $s_{21}=s_{32}=-a_{22}/2\le0$. The points associated to
$\SHP(A)$ are ${v^1}'=\ptop
0{\frac{a_{22}}{2}}{-2}=\ptop 2{2+\frac{a_{22}}{2}}0$ and
${v^3}'=\ptop{-2}{\frac{a_{22}}{2}}0$ and the slope of the segment $\overline{{v^3}'{v^1}'}$ is  $\frac{1}{2}$ (in $Z=0$),
independently of the precise value of $a_{22}$.
Then, any $a_{22}\ge0$ will do. We may take $a_{22}=0$ and we conclude that  $\CC$ is given by the
tropical polynomial $P=(-4)\odot X^{\odot2}\oplus
Y^{\odot2}\oplus Z^{\odot2}\oplus (-2)\odot X \odot
Y\oplus
 Y \odot Z\oplus  X \odot Z.$
\end{ex}

\section{Factorization of degree--two tropical polynomials}\label{sec:factorization}

A tropical polynomial $p$ (homogeneous or not) in any
number of variables is called \emph{reducible} if it is
the tropical product of two non--constant tropical
polynomials. A tropical hypersurface $\CC$ (affine or
projective) is called \emph{reducible} if it is the union
of  two hypersurfaces (affine or projective, accordingly)
$\CC_1,\CC_2$ with $\CC_1\neq\CC\neq\CC_2$. It is clear
that  the reducibility of a polynomial causes the
reducibility of the corresponding hypersurface but, in the
tropical setting, the converse is NOT true; see corollary
\ref{cor:irred_pol_red_con} below.

\m Let $P$ be  a homogeneous degree--two
tropical polynomial in three variables.  % We have $\sucesimp
%a\in\R$ and $\sucespar a\in\T$.
%and consider the associated
%matrix $A=D(A)\odot\SHP(A)\odot D(A)$.
The simplest example of a reducible polynomial arises when $a_{ij}=0$, for all $i\neq j$.
Then $P$ is the  square of the linear form
$\frac{a_{11}}{2}\odot X\oplus
\frac{a_{22}}{2}\odot Y\oplus\frac{a_{33}}{2}\odot Z,$ because in  tropical algebra the
\emph{freshman's dream}  $(a\oplus b)^{\odot n}=a^{\odot n}\oplus b^{\odot n}$ holds for all $n$!
The corresponding matrices and points are easy to compute:
$$A=\left(\begin{array}{ccc} a_{11}&\\ -\infty&a_{22}\\ -\infty&-\infty&a_{33}\\
\end{array}\right)=D^{\odot 2},$$
$\SHP(A)=I$ is the tropical identity matrix, $\SHP(A)^+$  is the zero matrix and the tropical conic $\CC(P)$ is a \emph{double line} with vertex at  $v=\frac{1}{2}[-a_{11},-a_{22},-a_{33}]$.

\begin{lem}\label{lem:fac}
$P$ is reducible if and only if $\SHP(P)$ is.
\end{lem}
\begin{proof} Consider the associated matrix $A=A(P)$.
The factorization $A=D\odot S\odot D$ corresponds to a change of variables
$[X,Y,Z]\mapsto [X',Y',Z']=[X,Y,Z]\odot D.$
\end{proof}

The former lemma allows us to reduce our discussion to the case $P=\SHP(P)$.

\begin{thm}
\label{thm:P=shpP} If  $P=\SHP(P)=\polshp
s$, then the  following statements hold.
\begin{enumerate}
\item If $-\infty\le s_{ij}<0$, for some
$i\neq j$, then $P$ is irreducible.
\item If $s_{ij}\ge0$, for all
$i\neq j$, then $P$ is reducible if and only if the maximum of
$\sucespar s$ equals the sum of the other two.
\end{enumerate}
\end{thm}

\begin{proof}
Up to tropical multiplication by a real
constant, a tropical factorization of $P$ must have the
form
\begin{equation}\label{eqn:factorization} \left(a\odot X\oplus b\odot Y\oplus
Z\right)\odot \left((-a)\odot X\oplus (-b)\odot Y\oplus
 Z\right),
\end{equation}
 for $a, b\in\R$, where
${s_{21}}=|a-b|$, ${s_{32}}=|b|$,
${s_{31}}=|a|\in\R_{\ge0}$. The irreducibility statement  now follows.
For the second statement, let us assume  that $s_{ij}\ge0$, for all
$i\neq j$ and, without loss of generality,
that $s_{31}=\max\{s_{21},s_{32},s_{31}\}$. Suppose that $s_{31}=s_{21}+s_{32}$.
Then we  take $a=s_{31}$ and $b={s_{32}}$, so that $P$ equals the product
(\ref{eqn:factorization}). The converse is easy.
\end{proof}

\begin{cor}\label{cor:irred_pol_red_con} If $P=\SHP(P)$ and $-\infty\le s_{ij}<0$, for
all $i\neq j$, then the polynomial $P$ is irreducible, but the conic $\CC(P)$  is a double
line.\qed
\end{cor}

%%By theorem \ref{lem:singular}, if $P=\polshp f$ with
%% $\sucespar f$ all non--negative and $P$ is reducible, then the matrix
%% $\SHP(A)=\SHP(A)^+$  is tropically singular. The converse
%% is NOT true.

\m Summing up, here is
a \label{proc:factor} \textbf{procedure to determine
whether  a given tropical degree--two homogeneous polynomial $P$ in three variables is reducible and,
in such a case, to
obtain a factorization}.
\begin{itemize}

\item Compute
the polynomial $\SHP(P)$ and decide whether it is reducible or not,
using theorem \ref{thm:P=shpP}.

\item If $\SHP(P)$ is reducible, we can factor it, as explained in the proof of theorem
\ref{thm:P=shpP}. Then,  a change of coordinates provides a factorization of $P$, by lemma \ref{lem:fac}.
\end{itemize}

\begin{ex}\label{ex:factorization} Let  $P=X^{\odot2}\oplus 12\odot
Y^{\odot2}\oplus Z^{\odot2}\oplus 7\odot X \odot Y\oplus
6\odot Y \odot Z\oplus 1\odot X \odot Z$.  The associated
matrices are
$$A=
\left(\begin{array}{ccc}
0\\
7&12\\
1&6&0\\
\end{array}\right),
\ D=
\left(\begin{array}{ccc}
0\\
-\infty&6\\
-\infty&-\infty&0\\
\end{array}\right),
\
S=S^+=
\left(\begin{array}{ccc} 0\\
1&0\\
1&0&0\\
\end{array}\right).$$ Then  $s_{21}=s_{31}=1$, $s_{32}=0$ and $\max\{\sucespar s\}=s_{31}=s_{21}+s_{32}$.
By theorem \ref{thm:P=shpP}, the polynomial
 $\SHP(P)=X^{\odot2}\oplus Y^{\odot2}\oplus Z^{\odot2}
 \oplus
1\odot X \odot Y\oplus Y \odot Z \oplus 1\odot X\odot Z$
is reducible and a factorization is given by
(\ref{eqn:factorization}) with $a=1$, $b=0$. Then the
translation given by the point  $[0,-6,0]$ provides the following factorization
$P=\left(1\odot X\oplus (-6)\odot Y\oplus Z\right)\odot
\left((-1)\odot X\oplus6\odot Y\oplus Z\right).$
\end{ex}

\end{document}